\documentclass[journal]{IEEEtran}

\ifCLASSINFOpdf
   \usepackage[pdftex]{graphicx}
\else
\fi
\usepackage{amsmath}
\usepackage{amssymb}
\usepackage{url}

\usepackage{algorithm,algorithmic}

\usepackage{array}
\usepackage{tabularx}

\usepackage{color}
\newcommand{\cm}{\color{black}}
\usepackage{cite} 

\begin{document}
%
\title{A Hierarchical Optimization Architecture for\\   Large-Scale Power Networks}
%
%
%

\author{Sungho Shin, Philip Hart, Thomas Jahns, and Victor M. Zavala

\thanks{S. Shin and V. M. Zavala are with the Department
of Chemical and Biological Engineering (e-mail: victor.zavala@wisc.edu) and P. Hart and T. Jahns are with the Department
of Electrical and Computer Engineering, University of Wisconsin-Madison, 1415 Engineering Dr, Madison, WI 53706, USA. }}
%
%

\markboth{}%
{Shell \MakeLowercase{\textit{et al.}}: Bare Demo of IEEEtran.cls for IEEE Journals}
%



\maketitle

\begin{abstract}
We present a hierarchical optimization architecture for large-scale power networks that overcomes limitations of fully centralized and fully decentralized architectures. The architecture leverages principles of multigrid computing schemes, which are widely used in the solution of partial differential equations on massively parallel computers. The top layer of the architecture uses a coarse representation of the entire network while the bottom layer is composed of a family of decentralized optimization agents each operating on a network subdomain at full resolution. We use an alternating direction method of multipliers (ADMM) framework to drive coordination of the decentralized agents. We show that state and dual information obtained from the top layer can be used to accelerate the coordination of the decentralized optimization agents and to recover optimality for the entire system. We demonstrate that the hierarchical architecture can be used to manage large collections of microgrids. 
\end{abstract}


%
\IEEEpeerreviewmaketitle

\section{Introduction}\label{sec:Intro}
\IEEEPARstart{P}{ower} networks are becoming increasingly difficult to manage due to the deployment of intermittent renewable power and the deployment of large collections of distributed energy resources. As complexity increases, the limitations of centralized optimization and control architectures are becoming increasingly evident. In particular, centralized architectures will not be capable of managing increasing amounts of sensor information and decisions. These limitations can be overcome by using decentralized architectures because this can mitigate communication and computation (decision-making) needs. However, decentralized architectures are limited in that coordination of large collections of agents can be slow, ultimately leading to robustness, stability, and economic performance issues. 

Hierarchical architectures provide a framework to overcome the challenges of centralized and decentralized decision-making. Such architectures have been recently explored in the context of model predictive control to handle behavior occurring at different timescales \cite{Zavala2016,Scattolini2007,Scattolini2009,Picasso2010,Rawlings2009}. Here, the fundamental idea is to use a supervisory layer that makes decisions over slow timescales and long horizons and a lower layer that conducts decisions over fast timescales and short horizons. Information in the form of state targets is used to ensure consistency between the control layers. Recently, we have also shown that hierarchical schemes can be used to make decisions over multiple {\em spatial} scales \cite{shin2018multi}. Such schemes are based on the observation that a hierarchical architecture shares similarities with multigrid computing schemes used in the solution of partial differential equations (PDEs) on massively parallel computers \cite{Hackbusch2013,Borzi2003}. In a spatial setting, the top layer uses a coarse representation of the entire spatial domain and the lower layer that is composed of a collection of decentralized agents, each operating on a spatial subdomain at full resolutions. State and dual information computed by the top layer guides and accelerate the coordination of the decentralized agents in the lower layer. The key insight, as noticed in the context of PDEs, is that decentralized agents are capable of mitigating local disturbances (also known as defects with high spatial frequencies) while the supervisory layer addresses global disturbances (with low spatial frequencies).  In the approach proposed in \cite{shin2018multi}, a Gauss-Seidel scheme is used to conduct coordination of the decentralized agents. Such an approach is intuitive and acts as a smoother of local disturbances but has limited convergence guarantees, particularly in nonconvex settings.  

{\cm In this work, we extend the hierarchical architecture proposed in \cite{shin2018multi} by (i) developing decomposition and coarsening strategy that can be applied to any type of graph structures, (ii) using an alternating direction of method of multipliers (ADMM) framework to perform the coordination, and (iii) demonstrating the framework with challenging large-scale AC optimal power flow (OPF) problems.} Our choice of ADMM is based on the observation that this is, in fact,  a Gauss-Seidel coordination scheme but applied to an augmented Lagrangian function of the system. {\cm ADMM has well-established convergence guarantees in a convex setting \cite{Tyrrell1976,Boyd2010} and in some restricted nonconvex settings \cite{wang2015global,hong2016convergence,Erseghe2014}}. We apply the hierarchical framework to power flow management problems over large power networks and argue that the proposed approach can coordinate large collections of distributed energy resources and microgrids. Under the proposed framework, decentralized agents operating over individual microgrids will seek to optimize local performance while coordinating with others to optimize system-wide performance. A supervisory layer ensures that fast coordination among the microgrids and network regions is achieved. We demonstrate that the ADMM framework performs satisfactorily in nonconvex network problems when used within a hierarchical architecture because this is aided by the supervisory layer.

The paper is organized as follows. In Section \ref{sec:setting} we present the notation and problem setting. In Section \ref{sec:ADMM} we present an overview of ADMM coordination. In Section \ref{sec:multi-grid} we formulate the supervisory layer and the hierarchical architecture. In Section \ref{sec:casestudy} we provide numerical experiments. 

\begin{figure}[!t]
\begin{center}
\includegraphics[width=.4\textwidth]{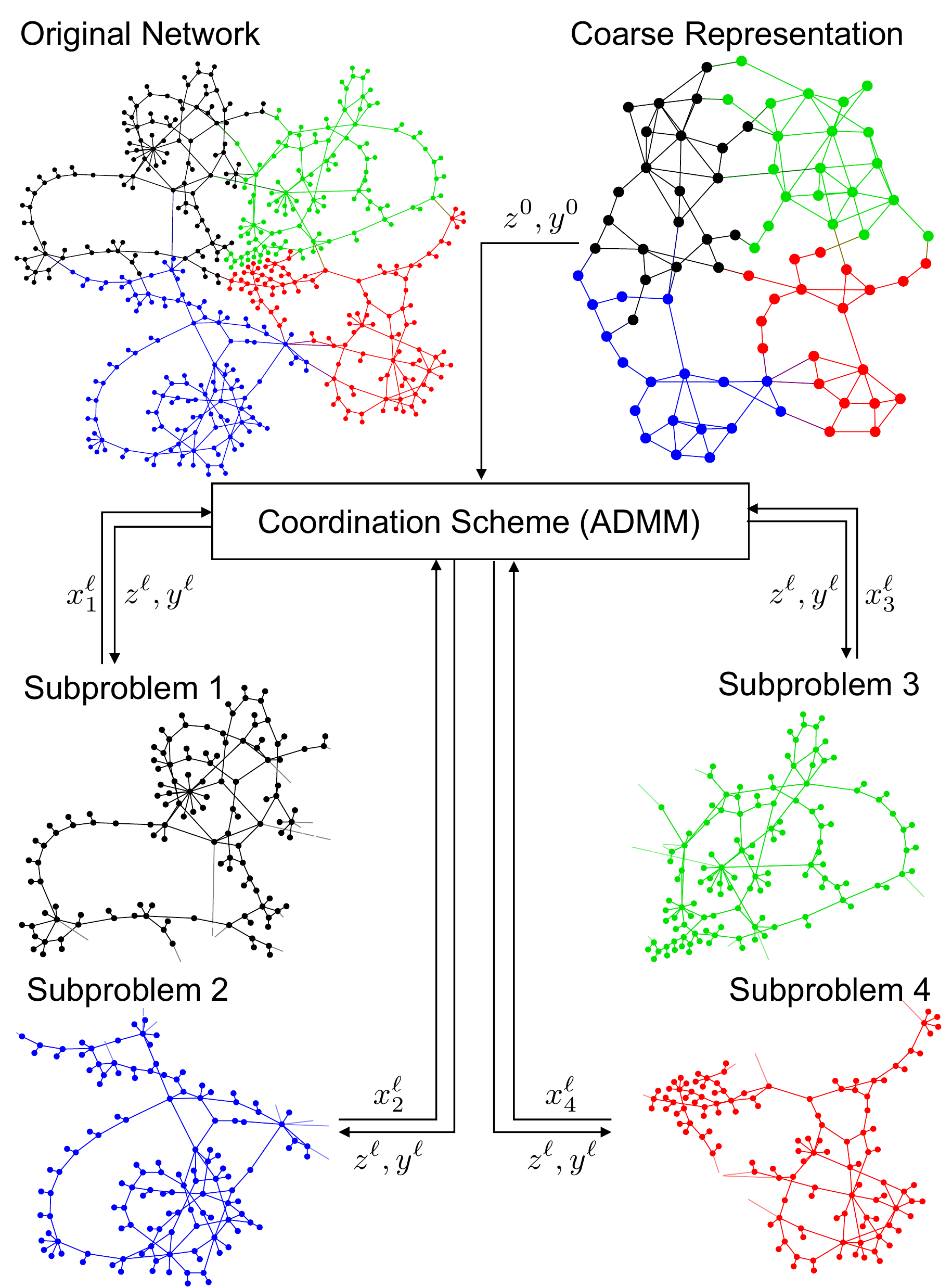}
\caption{Schematic of hierarchical architecture. A supervisory layer uses a coarse representation of the original network to guide coordination of decentralized agents that operate over individual partitions at full resolution. Network topology corresponds to that of instance case500\_tamu in the PGlib library \cite{pglib}.}
\label{fig:orig}
\end{center}
\end{figure}

\section{Problem Setting}\label{sec:setting}
\subsection{Basic Notation}
We consider all vectors as column vectors and use the syntax $(v_1,v_2,\cdots ,v_n) = \begin{bmatrix}v_1^T v_2^T \cdots v_n^T\end{bmatrix} ^T$. To represent vectors with irregular indices, we use the notation $\{x(i)\}_{i\in\mathcal{I}}:=(x(i_1),x(i_2),\cdots,x(i_n))$, where $\mathcal{I}:=\{i_1<i_2<\cdots<i_n\}$. Furthermore, we write $x_\mathcal{I}:=\{ x_i \}_{i\in\mathcal{I}}$ and $x(\mathcal{I}):=\{x(i)\}_{i\in\mathcal{I}}$. Euclidean vector norms are denoted as $\Vert\cdot\Vert$. The power networks discussed in this work are represented as undirected graphs $\mathcal{G}(\mathcal{V},\mathcal{E})$ where $\mathcal{V}$ is the set of nodes (buses) and $\mathcal{E}$ is the set of edges (lines). The edge between node $i$ and node $j$ is denoted by $\{i,j\}$. {Open and closed neighborhoods of node $i\in\mathcal{V}$ are denoted by $\mathcal{N} (i)$ and $\mathcal{N}[i]$, respectively, and defined as $\mathcal{N}(i):=\{j\in\mathcal{V}\,\mid \{i,j\}\in\mathcal{E}\}$ and $\mathcal{N}[i]:=\mathcal{N}(i)\cup\{i\}$, respectively.} When convenient, we use the notation $\mathcal{N}_{\mathcal{G}}$ to indicate that the neighborhoods are specific to graph $\mathcal{G}$. 

\subsection{Optimization Problem}
We consider the following network optimization problem:
\begin{subequations}\label{eqn:full}
\begin{align}
\min_{x}\quad & \sum_{i\in\mathcal{V}} f_{(i)} (x(i)) \\
\label{eqn:full_eq}\text{s.t.}\quad & g_{(i)}(x(\mathcal{N}[i])) =0 ,\quad i \in \mathcal{V}\\
\label{eqn:full_ineq} & h_{(i)}(x(\mathcal{N}[i]))\leq 0,\quad i \in \mathcal{V}
\end{align}
\end{subequations}
Here, the state variables $x(i)\in\mathbb{R}^{n_x(i)}$ are associated with node $i\in\mathcal{V}$. Similarly, the objective function $f_{(i)}$ and the constraint functions $g_{(i)},h_{(i)}$ are defined over each node  $i\in\mathcal{V}$. The notation $g_{(i)}(x(\mathcal{N}[i]))$ indicates that the constraint of node $i\in \mathcal{V}$ depends on its local states and on the states of all its neighbor nodes.  All functions are assumed to be at least twice continuously differentiable and potentially nonconvex. 

In the context of power networks, the states represent voltages, voltage angles, active and reactive power injections; the objective function is typically the total generation cost function or a tracking function that keeps state variables at desired levels; the equality constraint functions contain power flow equations; and the inequality constraint functions typically represent physical constraints such as power generation and voltage limits. Network data such as line admittances are embedded in the problem functions. The proposed framework is presented in the context of optimal power flow (in order to facilitate the presentation) but it is certainly not specific to this problem and can be used for other networks such as natural gas, water, and transportation networks as well as supply chains.

\subsection{Network Partitioning} 
We partition the entire node set of the original network into $K$ subsets of the form 
$\{\mathcal{V}_1,\mathcal{V}_2,\cdots,\mathcal{V}_K\}$. The set of partitions is defined as $\mathcal{K}:=\{1,2,\cdots,K\}$. We assume that the node subsets $\mathcal{V}_k$ are nonempty, disjoint, and $\mathcal{V}=\bigcup_{k\in\mathcal{K}} \mathcal{V}_k$.

The subnetworks associated with the partitions are coupled via {\em coupling nodes}. In order to distinguish such nodes, we define the following sets:
  \begin{align}\label{eqn:coupling}
    \overline{\mathcal{V}}_k &:=\bigcup_{k'\in\mathcal{K}\setminus\{k\}} \mathcal{N}[\mathcal{V}_k] \cap \mathcal{N}[\mathcal{V}_{k'}],\quad
    \overline{\mathcal{V}} := \bigcup_{k\in\mathcal{K}} \overline{\mathcal{V}}_k.
  \end{align}
Here, $\overline{\mathcal{V}}_k$ is the set of coupling nodes in $\mathcal{N}[\mathcal{V}_k]$, and $\overline{\mathcal{V}}$ is the set of all coupling nodes in the network. We use these sets to define local partition variables $x_k$ and global coupling variables $z$ associated with the coupling nodes:
\begin{subequations}
\begin{align}
  x_k&:=\{x_k(i)\}_{i\in\mathcal{N}[\mathcal{V}_k]},\quad k\in\mathcal{K}\\
  x_\mathcal{K}&:=(x_1,x_2,\cdots,x_K),\quad
  z:=\{z(i)\}_{i\in\overline{\mathcal{V}}}.
\end{align}
\end{subequations}
For simplicity, we write the full state vector as $x=x_\mathcal{K}$. 

The above definitions allow us to express \eqref{eqn:full} in the following {\em lifted form}:
\begin{subequations}\label{eqn:lifted}
\begin{align}
\min_{x,z}\quad &  \sum_{k\in\mathcal{K}}\sum_{i\in\mathcal{V}_k} f_{(i)} (x_k (i)) \\
\text{s.t.}\quad & (\lambda_k(i))\quad g_{(i)}(x_k(\mathcal{N}[i]))  =0 ,\quad i \in \mathcal{V}_k,\ k \in \mathcal{K} \label{eqn:lifted_eq}\\
 &(\nu_k(i))\quad h_{(i)}(x_k(\mathcal{N}[i]))\leq 0 ,\quad i \in \mathcal{V}_k,\ k \in \mathcal{K}\label{eqn:lifted_ineq}\\
&(y_k(i))\quad x_k(i)=z(i),\quad i \in \overline{\mathcal{V}}_k ,\ k \in \mathcal{K}.\label{eqn:lifted_link}
\end{align}
\end{subequations}
Here, the Lagrange multipliers associated with \eqref{eqn:lifted_eq}, \eqref{eqn:lifted_ineq}, and \eqref{eqn:lifted_link} are denoted by $\lambda_k(i)$, $\nu_k(i)$, and $y_k(i)$, respectively. By grouping variables by partitions, \eqref{eqn:lifted} can be written as:
\begin{subequations}\label{eqn:partition}
\begin{align}
\min_{x,z}\quad &  \sum_{k\in\mathcal{K}}f_k (x_k) \\
\text{s.t.}\quad & (\lambda_k)\quad g_k(x_k)  =0 ,\quad k \in \mathcal{K} \label{eqn:partition_eq}\\
 & (\nu_k)\quad h_k(x_k)\leq 0 ,\quad k \in \mathcal{K}\label{eqn:partition_ineq}\\
&(y_k)\quad A_k x_k+B_k z= 0,\quad k \in \mathcal{K}\label{eqn:partition_link},
\end{align}
\end{subequations}
where $A_k$ and $B_k$ are incidence matrices. The linking constraints \eqref{eqn:partition_link} are defined for the coupling nodes $\overline{\mathcal{V}}_k$ of each partition $k\in\mathcal{K}$.
To achieve compact notation, we define:
\begin{subequations}\label{eqn:cpts}
  \begin{align}
    g(x) &=(g_1(x),\cdots,g_K(x)),\; h(x) =(h_1(x),\cdots,h_K(x))\\
    A&=\mathop{\textrm{blkdiag}}(A_1 ,\cdots ,A_K),\; B=\left[B_1^T\; \cdots\; B_K^T\right]^T.
  \end{align}
\end{subequations}
The  incidence matrices $A$ and $B$ have full row rank and full column rank, respectively. By defining $\mathcal{S}:=\{x \mid  g(x) =0 \text{ and } h(x)\leq 0\}$, we obtain the following compact form: 
\begin{subequations}\label{eqn:lifted_cpt}
\begin{align}
\min_{x\in\mathcal{S},z}\quad & f(x)\label{eqn:lifted_cpt_obj} \\
\text{s.t.}\quad &(y)\; Ax+Bz=0.\label{eqn:lifted_cpt_link}
\end{align}
\end{subequations}
Note that \eqref{eqn:lifted}, \eqref{eqn:partition}, and \eqref{eqn:lifted_cpt} are different representations of the same problem. We denote the Lagrange multiplier associated with $g(x)=0$ and $h(x)\leq 0$ as $\lambda$ and $\nu$, respectively.

\section{Coordination Scheme}\label{sec:ADMM}
Problem \eqref{eqn:lifted_cpt} can be solved in a decentralized manner by using an ADMM scheme. The convergence properties of ADMM have been studied extensively in convex optimization settings, but convergence properties in nonconvex settings are still not well understood and require experimental testing. {\cm Recent studies have demonstrated, however, that ADMM works well in nonconvex problems because its performance resembles that of an augmented Lagrangian method\cite{Erseghe2014,sun2013fully,rodriguez2018benchmarking,guo2017case}.}

\subsection{Overview of ADMM}
The {\em partial} augmented Lagrangian of \eqref{eqn:lifted_cpt} is given by \cite{Bertsekas1982}:
\begin{align}\label{eqn:augLag}
\mathcal{L}_\rho (x,z,y):=f(x)+y^T(Ax+Bz)+\frac{1}{2}\rho\Vert Ax+Bz \Vert^2,
\end{align}
where $y$ is the Lagrange multiplier of \eqref{eqn:lifted_cpt_link} and $\rho\geq 0$ is a penalty parameter. Note that only the coupling constraints \eqref{eqn:lifted_cpt_link} are incorporated into the partial augmented Lagrangian. In other words, the inner (potentially nonconvex) constraints $x\in\mathcal{S}$ are not incorporated. The ADMM scheme solves the following problems over $\ell=0,1,2,\cdots$ coordination steps:
\begin{subequations}\label{eqn:ADMM_alg}
\begin{align}
x^{\ell+1}&=\mathop{\textrm{argmin}}_{x\in\mathcal{S}} \mathcal{L}_\rho (x,z^\ell,y^\ell) \label{eqn:ADMM_x}\\
z^{\ell+1}&=\mathop{\textrm{argmin}}_{z} \mathcal{L}_\rho (x^{\ell+1},z,y^\ell) \label{eqn:ADMM_z}\\
y^{\ell+1}&=y^\ell+\rho(Ax^{\ell+1}+Bz^{\ell+1}) \label{eqn:ADMM_y}
\end{align}
\end{subequations}

Problem \eqref{eqn:ADMM_x} can be fully split into individual network partitions, as it uses fixed values of the coupling states $z^\ell$ and duals $y^\ell$. Accordingly, the splitted individual problems can be solved in parallel for each $k\in\mathcal{K}$. Using the notation of \eqref{eqn:partition}, the subproblem over partition $k$ (corresponding to step \eqref{eqn:ADMM_x}) can be written as:
\begin{subequations}\label{eqn:ADMM_local}
\begin{align}
\min_{x_k} \quad& f_k(x_k)+(y_k^\ell)^T\left(A_k x_k+B_k z^\ell\right)+\frac{1}{2}\rho \Vert A_k x_k+B_k z^\ell \Vert^2 \\
\text{s.t.}\quad & g_k(x_k) =0  \label{eqn:ADMM_local_eq} \\
 & h_k(x_k)\leq 0.\label{eqn:ADMM_local_ineq}
\end{align}
\end{subequations}
This structure reveals that $z^\ell$ act as {\em state targets} that each partition seeks to follow and that $y^\ell$ act as {\em prices for information} on the coupling variables.

Coordination of the network partitions is achieved by updating the coupling states $z$ with \eqref{eqn:ADMM_z} and dual $y$ with \eqref{eqn:ADMM_y}. Subproblem \eqref{eqn:ADMM_z} is an unconstrained convex quadratic program and has a closed-form solution of the form:
\begin{align}\label{eqn:z_closed}
z^{\ell+1}=-(B^TB)^{-1}B^T( A x^{\ell+1} + \frac{1}{\rho}y^\ell),\; \ell=0,1\cdots.
\end{align}
{The nonsingularity of $B^T B$ comes from the fact that $B$ has full column rank. From \eqref{eqn:By=0} in Appendix \ref{apx:KKT} we have that $B^T y^\ell=0$ for $\ell=1,2,\cdots$. Thus we have:
\begin{align}
z^{\ell+1}=-(B^TB)^{-1}B^T A x^{\ell+1},\; \ell=1,2,\cdots.
\end{align}}
Using the notation of \eqref{eqn:lifted}, one can show that this reduces to:
\begin{align}\label{eqn:ADMM_global}
z^{\ell+1}(i)=\frac{1}{|\mathcal{K}_i|}\sum_{k\in\mathcal{K}_i} x_k^{\ell+1}(i),\; \forall i \in \overline{\mathcal{V}},\;\ell=1,2,\cdots.
\end{align}
and where $\mathcal{K}_i:=\{k\in\mathcal{K}\mid i\in\overline{\mathcal{V}}_k\}$. This reveals that the coupling states are updated (coordinated) by {\em averaging} the states associated with the partitions connected to the coupling node $i$ (given by the set $\mathcal{K}_i$). Step \eqref{eqn:ADMM_y} can be written as follows.
\begin{equation}\label{eqn:ADMM_dual}
y_k^{\ell+1}(i)=y_k^{\ell}(i)+\rho (x_k^{\ell+1}(i)-z^{\ell+1}(i)),\; \forall i \in \overline{\mathcal{V}}_k,k\in\mathcal{K}
\end{equation}
This reveals that this coordination step seeks to close the gap between states across partitions.  The structure of the state and dual updates indicates that information only needs to be shared between neighboring nodes. 

\subsection{Monitoring Convergence}\label{subsec:opt_cond}
We now proceed to derive conditions that allow us to determine if optimality has been achieved by the ADMM scheme. The Lagrangian for \eqref{eqn:lifted_cpt} is: 
\begin{align}
\mathcal{L} (x,z,y,\lambda,\nu):=&f(x)+{\lambda}^Tg(x)+{\nu}^Th(x)+y^T( A x+B z).
\end{align}
{\cm It is well-known that, in a convex setting, the primal and dual residuals of the first-order optimality conditions of  \eqref{eqn:lifted_cpt} at coordination step $\ell=1,2,\cdots$ are given by \cite{Boyd2010}:
\begin{subequations}\label{eqn:res}
\begin{align}
r^\ell &:= Ax^\ell+Bz^\ell\label{eqn:r}\\
s^\ell &:= \rho A^T B(z^{\ell} - z^{\ell-1}) \label{eqn:s}.
\end{align}
\end{subequations}
Using the notation of \eqref{eqn:lifted}, these residuals can be written as:
\begin{subequations}
\begin{align}
r^\ell_k(i)&=x_k^\ell(i)-z^\ell(i), \quad i\in\overline{\mathcal{V}}_k,\ k\in\mathcal{K}\\
s_k^\ell(i)&=\begin{cases}
\rho\left(-z^{\ell}(i)+z^{\ell-1}(i)\right)&\text{if } i\in\overline{\mathcal{V}}_k
\\
0&\text{if }i\in\mathcal{N}[\mathcal{V}_k]\setminus\overline{\mathcal{V}}_k.
\end{cases}
\end{align}
\end{subequations}

Now consider $\ell$-th coordination step $(z^\ell,y^\ell)$ of ADMM \eqref{eqn:ADMM_alg}; let $(x^{\ell+1},\lambda^{\ell+1},\nu^{\ell+1})$ be the primal-dual solution of \eqref{eqn:ADMM_x} such that satisfies some constraint qualification, $z^{\ell+1}$ be the solution of \eqref{eqn:ADMM_z}, and $y^{\ell+1}$ be obtained from \eqref{eqn:ADMM_y}. In Appendix \ref{apx:KKT}, we show that the standard form of the residuals \eqref{eqn:res} can also be used in a partial augmented Lagrangian setting \eqref{eqn:augLag} with a nonconvex feasible set.} In particular, the residual of the first order optimality conditions of \eqref{eqn:lifted_cpt} evaluated at $(x^{\ell+1},z^{\ell+1},y^{\ell+1},\lambda^{\ell+1},\nu^{\ell+1})$ can be expressed by:
\begin{subequations}\label{eqn:residual}
\begin{align}
\nabla_x \mathcal{L} (x^{\ell+1},z^{\ell+1},y^{\ell+1},\lambda^{\ell+1},\nu^{\ell+1}) &=  s^{\ell+1} \label{eqn:stationary_1} \\
\nabla_z \mathcal{L} (x^{\ell+1},z^{\ell+1},y^{\ell+1},\lambda^{\ell+1},\nu^{\ell+1}) &= 0  \label{eqn:stationary_2}\\
g(x^{\ell+1}) &= 0 \label{eqn:pr_fs_1} \\
h(x^{\ell+1}) &\leq 0 \label{eqn:pr_fs_2}\\
Ax^{\ell+1}+Bz^{\ell+1} &= r^{\ell+1} \label{eqn:pr_fs_3}\\
\nu^{\ell+1}&\geq 0\label{eqn:du_fs}\\
\text{diag}(\nu^{\ell+1}) h(x^{\ell+1}) &= 0 \label{eqn:cmpl_slk} 
\end{align}
\end{subequations}
From these expressions it becomes clear that if $\Vert r^\ell \Vert, \Vert s^\ell \Vert \rightarrow 0$, then a stationary point of \eqref{eqn:lifted_cpt} is obtained. Consequently, all that is needed to check convergence are the primal and dual residuals $r^\ell$ and $s^\ell$, which are easy to evaluate from \eqref{eqn:res} by sharing information among the neighboring nodes. {\cm It is important to note that, in nonconvex settings, the conditions $\Vert r^\ell\Vert=\Vert s^\ell\Vert =0$ only guarantee stationarity.}

As with most iterative algorithms, the ADMM coordination scheme approaches the solution asymptotically. Consequently, it is expected that a good initial estimate for the coupling states and duals $(z^0,y^0)$ will aid the algorithm to reach optimality in fewer coordination steps \cite{rodriguez2018benchmarking}.  

\section{Hierarchical Architecture}\label{sec:multi-grid} 
A hierarchical architecture can be used to guide and accelerate the ADMM coordination scheme. This is done by obtaining estimates of the state and dual variables using a {\em coarse and tractable} representation of the network. 

\subsection{Partition/Subpartition Identification}

To generate a coarse representation for the supervisory layer, we construct network {\em subpartitions} and aggregate the nodes in each subpartition. We use the term subpartition to emphasize that each subpartition is a subset of some partition $\mathcal{V}_k$. The aggregation procedure aims to drastically reduce the dimensionality of the problem to enable tractability of the supervisory layer.
To identify the partition/subpartition structure, we use a generic graph partitioning method based on multilevel $k$-way partitioning \cite{karypis2000multilevel}. Physics-based graph partitioning methods, for example, {\em coherency or similarity} identification method \cite{Demarco1995}, can also be used.

\begin{figure*}
  \centering
  \includegraphics[width=\textwidth]{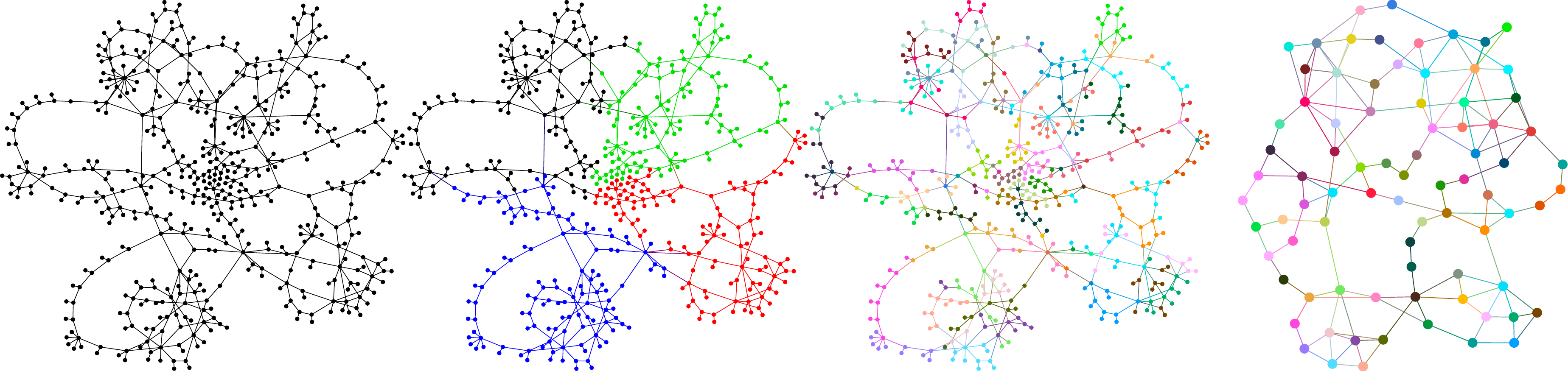}
  \caption{Left to right: Original (unpartitioned) network, partitioned network, subpartitions, coarse network (based on instance case500\_tamu in PGlib \cite{pglib}).}
  \label{fig:partitioning}
\end{figure*}

The graph partitioning is applied to first obtain the partition structure $\{V_1,V_2,\cdots,V_K\}$, and applied to each partition again to obtain the subpartition structure $\{\tilde{V}_1,\cdots,\tilde{V}_{K^c}\}$. The partition/subpartition identification and coarsening procedure are illustrated in Fig. \ref{fig:partitioning}. The coarse network is formed by aggregating the nodes in the same subpartition into a single node and aggregating the edges that connect the same pair of subpartitions into a single edge. The resulting coarse network misses the high-resolution details of the network but is able to capture the low-resolution, coarse behavior of the network. By formulating a coarse problem associated with the coarse network, we can construct a problem that can provide an approximate solution of the original problem.

\subsection{Coarse Graph}
The following discussion describes a scheme to aggregate nodes and edges to construct the coarse node set $\mathcal{V}^c$ and coarse edge set $\mathcal{E}^c$. With the aggregated node and edge sets, we aim to define the coarse graph $\mathcal{G}^c(\mathcal{V}^c,\mathcal{E}^c)$. 

Let $\{\tilde{\mathcal{V}}_1,\tilde{\mathcal{V}}_2,\cdots,\tilde{\mathcal{V}}_{K^c}\}$ be a collection of nonempty, disjoint subsets of $\mathcal{V}$ such that $\mathcal{V}=\bigcup_{i^c\in\mathcal{V}^c}\tilde{\mathcal{V}}_{i^c}$ holds, where $\mathcal{V}^c:=\{1,2,\cdots,K^c\}$. We assume that $\{\tilde{\mathcal{V}}_1,\tilde{\mathcal{V}}_2,\cdots,\tilde{\mathcal{V}}_{K^c}\}$ forms a subpartition structure of $\mathcal{V}$. That is,
\begin{align}\label{eqn:subparinpar}
\forall i^c\in\mathcal{V}^c, \quad \exists k\in\mathcal{K} \quad \text{s.t.}\quad \hat{\mathcal{V}}_{i^c} \subseteq \mathcal{V}_k.
\end{align}
We identify the subpartition index set $\mathcal{V}^c$ as the set of aggregated nodes. The set of aggregated edges $\mathcal{E}^c$ can be defined as:
\begin{align}
\mathcal{E}^c:=\{\{i^c, j^c\}\mid \exists \{i,j\}\in\mathcal{E}:i\in\tilde{\mathcal{V}}_{i^c},\; j\in\tilde{\mathcal{V}}_{j^c},\; i^c\neq j^c\}.
\end{align}
Now we consider $\mathcal{V}^c$ and $\mathcal{E}^c$ as the node and edge sets of the coarse graph $\mathcal{G}^c$. Thus, the coarse graph $\mathcal{G}^c(\mathcal{V}^c,\mathcal{E}^c)$ is constructed.

To construct the hierarchical architecture, it is necessary to {\em communicate information} between the layers. To do so, we create a mapping from the fine to the coarse space $\varphi:\mathcal{V}\rightarrow \mathcal{V}^c$. In particular, we consider a mapping of the form:
\begin{align}
\varphi(i)=i^c,\quad \text{if}\quad i\in {\mathcal{V}}_{i^c}
\end{align}
for $i\in\mathcal{V}$. The mapping $\varphi(\cdot)$ is well-defined because ${\tilde{\mathcal{V}}}_1,{\tilde{\mathcal{V}}}_2,\cdots,\tilde{\mathcal{V}}_{K^c}$ are disjoint and their union is $\mathcal{V}$.

Using the original partition structure $\{\mathcal{V}_1,\cdots,\mathcal{V}_K\}$, one can partition the coarse graph. 
The partitioned node sets $\mathcal{V}^c_k$, coupling node set $\overline{\mathcal{V}}^c_k$, and overall coupling node set $\overline{\mathcal{V}}^c$ for the coarse graph can be defined as follows.
\begin{subequations}
  \begin{align}
    \mathcal{V}^c_k &:= \{i^c\in\mathcal{V}^c \mid \hat{\mathcal{V}}_{i^c}\subseteq \mathcal{V}_k\}\\
    \overline{\mathcal{V}}^c_k &:=\bigcup_{k'\in\mathcal{K}\setminus\{k\}} \mathcal{N}_{\mathcal{G}^c}[\mathcal{V}^c_k] \cap \mathcal{N}_{\mathcal{G}^c}[\mathcal{V}_{k'}^c],\;
\overline{\mathcal{V}}^c := \bigcup_{k\in\mathcal{K}} \overline{\mathcal{V}}^c_k.
  \end{align}
\end{subequations}

\subsection{Coarse Problem Formulation}
 We express the coarse representation of \eqref{eqn:lifted} as:
\begin{subequations}\label{eqn:coarse}
\begin{align}
\min_{x^c,z^c}&\;  \sum_{k\in\mathcal{K}}\sum_{i^c \in\mathcal{V}^c_k} f_{c,(i^c)} (x^c_k (i^c)) \\
\text{s.t.}\; &(\lambda^c_k(i^c))\; g_{c,(i^c)}  (x^c_k(\mathcal{N}_{\mathcal{G}^c}[i^c])) =0,\; i^c \in \mathcal{V}^c_k,\ k \in \mathcal{K} \label{eqn:coarse_eq}\\
 &(\nu ^c_k(i^c))\; h_{c,(i^c)}(x^c_k(\mathcal{N}_{\mathcal{G}^c}[i^c]))\leq 0 ,\; i^c \in \mathcal{V}^c_k,\ k \in \mathcal{K}\label{eqn:coarse_ineq}\\
 & (y^c_k(i^c))\; x^c_k(i^c)-z^c(i^c) =0 ,\; i^c \in \overline{\mathcal{V}}^c_k ,\ k \in \mathcal{K},\label{eqn:coarse_link}
\end{align}
\end{subequations}
where $x^c$ denote the coarse state variables and $z^c$ denote the coarse coupling variables. The coarse objective and constraint functions are defined as:
\begin{subequations}\label{eqn:coarse_functions}
  \begin{align}
    f_{c,(i^c)}(x^c_k(i^c)) &= \sum_{i \in \tilde{\mathcal{V}}_{i^c}}f_{(i)}(x(i))\label{eqn:crsf}\\
    g_{c,(i^c)}(x^c_k(\mathcal{N}_{\mathcal{G}^c}[i^c])) &= \sum_{i \in \tilde{\mathcal{V}}_{i^c}}g_{(i)}(x(\mathcal{N}_{\mathcal{G}}[i]))\label{eqn:crsg}\\
    h_{c,(i^c)}(x^c_k(\mathcal{N}_{\mathcal{G}^c}[i^c])) &= \sum_{i \in \tilde{\mathcal{V}}_{i^c}}h_{(i)}(x(\mathcal{N}_{\mathcal{G}}[i]))\label{eqn:crsh}
\end{align}
\end{subequations}
for  $i^c \in \mathcal{V}^c_k,\ k \in \mathcal{K}$, where we enforce $x(i) = x^c(\varphi(i))$ for each $i\in\mathcal{V}$.

The primal-dual solution $(x^c,z^c,y^c,\lambda^c,\nu^c)$ of the coarse problem \eqref{eqn:coarse} can be mapped back to the fine space to initialize the ADMM procedure. The mapping is:
  \begin{subequations}\label{eqn:projection}
    \begin{align}
      x_k^0 (i)&= x_k^c (\varphi(i)),\quad k\in\mathcal{K},\; i\in\mathcal{V}_k\\
      z^0 (i)&= z^c (\varphi(i)),\quad i\in\overline{\mathcal{V}}\\
      y_k^0 (i)&= y_k^c (\varphi(i)),\quad k\in\mathcal{K},\; i\in\overline{\mathcal{V}}_k.
    \end{align}
\end{subequations}
Note that $x^0$ can be used to warm start the solution of \eqref{eqn:ADMM_x} and $z^0$ and $y^0$ can be used as a starting point of ADMM procedure \eqref{eqn:ADMM_alg}. This can be interpreted that the state and dual information of the coarse problem is used to guide the decentralized ADMM coordination scheme. 

The structure of \eqref{eqn:coarse} reveals that the coarse problem can be solved in a centralized manner (all at once as a general optimization problem) or by using a decentralized scheme such as ADMM (when centralized solution of the coarse problem is not possible). Consequently, the proposed hierarchical architecture can embed multiple layers of increasing coarseness as we move up the hierarchy.

{\cm We highlight that, in many practical cases, the dimensions of $g_{(i)}(\cdot)$ and $h_{(i)}(\cdot)$ may not match. Thus, the summation in \eqref{eqn:crsg}-\eqref{eqn:crsh} may not be well-defined. In such a case, domain-specific knowledge needs to be applied. In particular, with the standard OPF formulation, there exist inequality constraints associated with edges, and they make the dimensions of $h_{(i)}$ non-uniform. In such a case, constraints on vertices and constraints on edges can be handled separately. For example, suppose that the inequality constraints take the following form:
\begin{align}
  h(x)\leq 0 \iff 
  \begin{cases}
    h^v_{(i)} (x(i)) \leq 0&\forall i\in \mathcal{V}\\
    h^e_{\{i,j\}}(x(i),x(j))) \leq 0&\forall \{i,j\}\in \mathcal{E},
  \end{cases}
\end{align}
and the dimensions of $h^v_{(i)}(\cdot)$ and $h^e_{\{i,j\}}(\cdot)$ are uniform throughout $i\in\mathcal{V}$ and $\{i,j\}\in\mathcal{E}$, respectively.
We may consider the following aggregation scheme:
\begin{subequations}
  \begin{align}
    h^v_{c,(i^c)}(x^c(i^c)) &= \sum_{i \in {\mathcal{V}}_{i^c}}h^v_{(i)}(x(i)))\label{eqn:crsh-v}\\
    h^e_{c,\{i^c,j^c\}}(x^c(i^c),x^c(j^c)) &=
    \sum_{\{i,j\} \in {\mathcal{E}}_{i^c,j^c}}h^e_{\{i,j\}}(x(i),x(j))\label{eqn:crsh-e}
  \end{align}
\end{subequations}
where we enforce $x(i) = x^c(\varphi(i))$ for each $i\in\mathcal{V}$ and
\begin{align}\label{eqn:Ejc}
  \mathcal{E}_{i^c,j^c}:=\{\{i,j\}\in\mathcal{E}\mid i \in\tilde{\mathcal{V}}_{i^c} \text{ and } j\in\tilde{\mathcal{V}}_{j^c}\}.
\end{align}
Finally, one can use the following coarse constraints:
\begin{subequations}
  \begin{align}
    & h^v_{c,(i^c)}(x^c(i^c))\leq 0 , \quad \forall i^c \in \mathcal{V}^c\\
    &h^e_{c,\{i^c,j^c\}}(x^c(i^c),x^c(j^c))\leq 0 , \quad\forall\{i^c,j^c\} \in \mathcal{E}^c.
  \end{align}
\end{subequations}
A detailed derivation of coarse OPF is given in Appendix \ref{apx:agg}.}

\subsection{Hierarchical Optimization Scheme}
The hierarchical scheme is summarized in Algorithm \ref{alg:full}. We highlight that the local states for the partitions are computed in a fully decentralized manner. The solutions are then used to update the coupling states and adjoint variables.  The fact that the coarse problem can also be solved by using an ADMM scheme reveals that the proposed hierarchical architecture can be generalized to contain multiple levels (of increasing coarseness as we move up the hierarchy). This provides a mechanism to handle extremely large networks. In this work we limit the discussion to two levels.  We also note that the coarse level can update the targets of the fine ADMM layer at different times (e.g., when the primal and dual residuals are too large to be handled by ADMM efficiently). In other words, one can determine suitable threshold values for the primal and dual residuals that trigger an update from the coarse layer. In this work we only consider the case in which the coarse level provides targets at the first coordination step (i.e., $\ell=0$).
\begin{algorithm}[!htpb]
 \caption{Hierarchical optimization scheme.}
 \label{alg:full}
\begin{algorithmic}
 \STATE $(z^c,y^c)\leftarrow$ Solve coarse problem \eqref{eqn:coarse}
  \STATE $\ell\leftarrow 0,r^\ell\leftarrow \infty,s^\ell\leftarrow \infty$ Initialize ADMM coordination.
 \STATE $(z^\ell,y^\ell)\leftarrow$ Map coarse solution to fine space using \eqref{eqn:projection}
\\ 
 \WHILE{$\Vert r^\ell \Vert \geq \epsilon^{\text{pr}}\text{ and }\Vert s^\ell \Vert \geq \epsilon^{\text{du}}$}
\FOR{{\cm (in parallel)} $k\in\mathcal{K}$}
\STATE $x_k^{\ell+1}\leftarrow$ Solve subproblem \eqref{eqn:ADMM_local} with $(x^\ell,z^\ell,y^\ell)$.
\ENDFOR
\STATE $z^{\ell+1}\leftarrow$ Update coupling states using \eqref{eqn:ADMM_global}.
\STATE $y^{\ell+1}\leftarrow$ Update dual variables using \eqref{eqn:ADMM_dual}.
\STATE $r^{\ell+1}\leftarrow$ Compute primal feasibility using \eqref{eqn:r}.
\STATE $s^{\ell+1}\leftarrow$ Compute dual feasibility using \eqref{eqn:s}.
\STATE $\ell\leftarrow \ell+1$
\ENDWHILE
 \end{algorithmic}
\end{algorithm}

\section{Case Studies}\label{sec:casestudy}%
In this section, we demonstrate that the proposed hierarchical architecture works for optimal power flow problems for different large-scale networks. ADMM was first applied to OPF problems in \cite{Kim1997}, and many recent works have investigated various modifications \cite{Erseghe2014,Magnusson2015,Phan2014,Peng2014}. These works show that, while decentralized ADMM can tackle complex problems, a large number of coordination steps are often required. Our case study seeks to demonstrate that the supervisory layer can be used to accelerate the convergence of ADMM. All results can be reproduced using the code provided in \url{https://github.com/zavalab/JuliaBox/HierarchicalADMM}.

{\cm We consider a standard OPF problem \cite{coffrin2016qc}:
  \begingroup
  \allowdisplaybreaks
  \begin{subequations}\label{eqn:cs}
    \begin{align}
      \label{eqn:cs_a}\min_{{V,\theta,P,Q}}\ &\sum_{j\in\mathcal{W}}  c_{j} P(j)\\
      \label{eqn:cs_b}\text{s.t.}\
      &P_L(i)+\sum_{j\in \mathcal{W}_i} P(j)=\\
      &\sum_{j\in \mathcal{N}[i]}G_{ij}V (i)V (j)\cos\left(\theta (i)-\theta (j)\right)\nonumber\\
      &+\sum_{j\in \mathcal{N}[i]}B_{ij}V (i)V (j)\sin\left(\theta (i)-\theta (j)\right),\; i\in\mathcal{V}\nonumber\\
      \label{eqn:cs_c}
      &Q_L(i)+\sum_{j\in \mathcal{W}_i} Q(j)=\\
      &\sum_{j \in \mathcal{N}[i]}G_{ij}V (i)V (j)\sin\left(\theta (i) - \theta (j)\right)\nonumber\\
      &-\sum_{j \in \mathcal{N}[i]}B_{ij}V (i)V (j)\cos\left(\theta (i) - \theta (j)\right),\; i\in\mathcal{V}\nonumber\\
      \label{eqn:cs_d}& P^{\min}(j) \leq P(j)\leq P^{\max}(j),\;j\in \mathcal{W}\\
      \label{eqn:cs_e}& Q^{\min}(j) \leq Q(j)\leq Q^{\max}(j),\;j\in \mathcal{W}\\
      \label{eqn:cs_f}& \theta^{\min} \leq \theta(i)-\theta(j) \leq \theta^{\max} ,\; \{i,j\}\in \mathcal{E}\\
      \label{eqn:cs_g}& V^{\min}\leq V(i) \leq V^{\max},\; i\in \mathcal{V}\\
      \label{eqn:cs_h}& \theta(i)=0,\; i\in \mathcal{V}^{\text{ref}}.
    \end{align}
  \end{subequations}
  \endgroup
  The set of generators are represented by $\mathcal{W}$; the set of generators that are connected to node $i\in\mathcal{V}$ is represented by $\mathcal{W}_i$ (note that $\mathcal{W}=\bigcup_{i\in\mathcal{V}}\mathcal{W}_i$); the set of reference nodes are denoted by $\mathcal{V}^{\text{ref}}$. The voltage amplitude and voltage angle of node $i\in\mathcal{V}$ are denoted by $V(i)$ and $\theta(i)$, respectively; active power and reactive power generated by generator $j\in\mathcal{W}$ are denoted by $P(j)$ and $Q(j)$, respectively; the real part (conductance) and the imaginary part (susceptance) of $(i,j)$th component of admittance matrix is denoted by $G_{ij}$ and $B_{ij}$, respectively. The unit generation cost of generator $j\in\mathcal{W}$ is denoted by $c_{j}$; the active and reactive power load of node $i\in\mathcal{V}$ is denoted by $P_L(i)$ and $Q_L(i)$, respectively; the maximum and minimum voltage angle separation and voltage amplitude are denoted by $\theta_{\max}$, $\theta_{\min}$, $V_{\max}$, and $V_{\min}$, respectively; the maximum and minimum active and reactive power generation at generator $j\in\mathcal{W}$ are represented by $P^{\max}(j)$, $P^{\min}(j)$, $Q^{\max}(j)$, and $Q^{\min}(j)$, respectively. Here, the state variable of node $i\in\mathcal{V}$ is defined as $x(i)=\left(V(i),\theta(i),P(\mathcal{W}_i),Q(\mathcal{W}_i)\right)$.

  The objective function \eqref{eqn:cs_a} considers the total generation cost over the whole network. The standard polar form AC power flow equations \eqref{eqn:cs_b}-\eqref{eqn:cs_c} are used. Constraints \eqref{eqn:cs_d}-\eqref{eqn:cs_g} represent physical and stability limits. They include limits for active and reactive power generation \eqref{eqn:cs_d}-\eqref{eqn:cs_e}, voltage angle separations \eqref{eqn:cs_f}, and voltage amplitudes \eqref{eqn:cs_g}. Other types of constraints, such as maximum and minimum power flow constraints, are neglected for simplicity. The voltage angles of the reference nodes are fixed with \eqref{eqn:cs_h}. To make the problem always feasible, an artificial slack generator with a sufficiently large generation cost is placed in each node.}

\begin{figure*}[!htbp]
  \begin{center}
    \includegraphics[width=.19\textwidth]{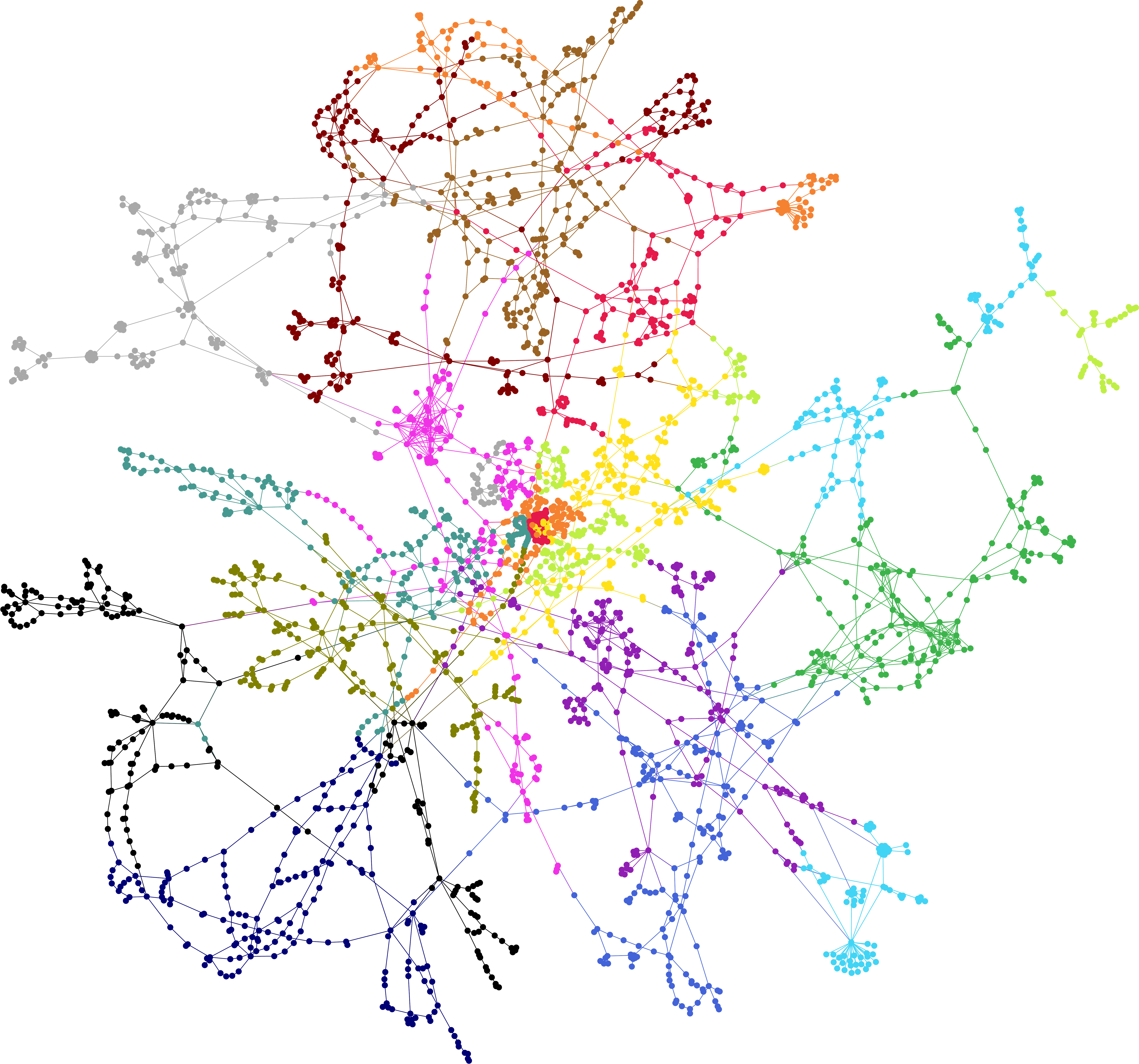}
    \includegraphics[width=.19\textwidth]{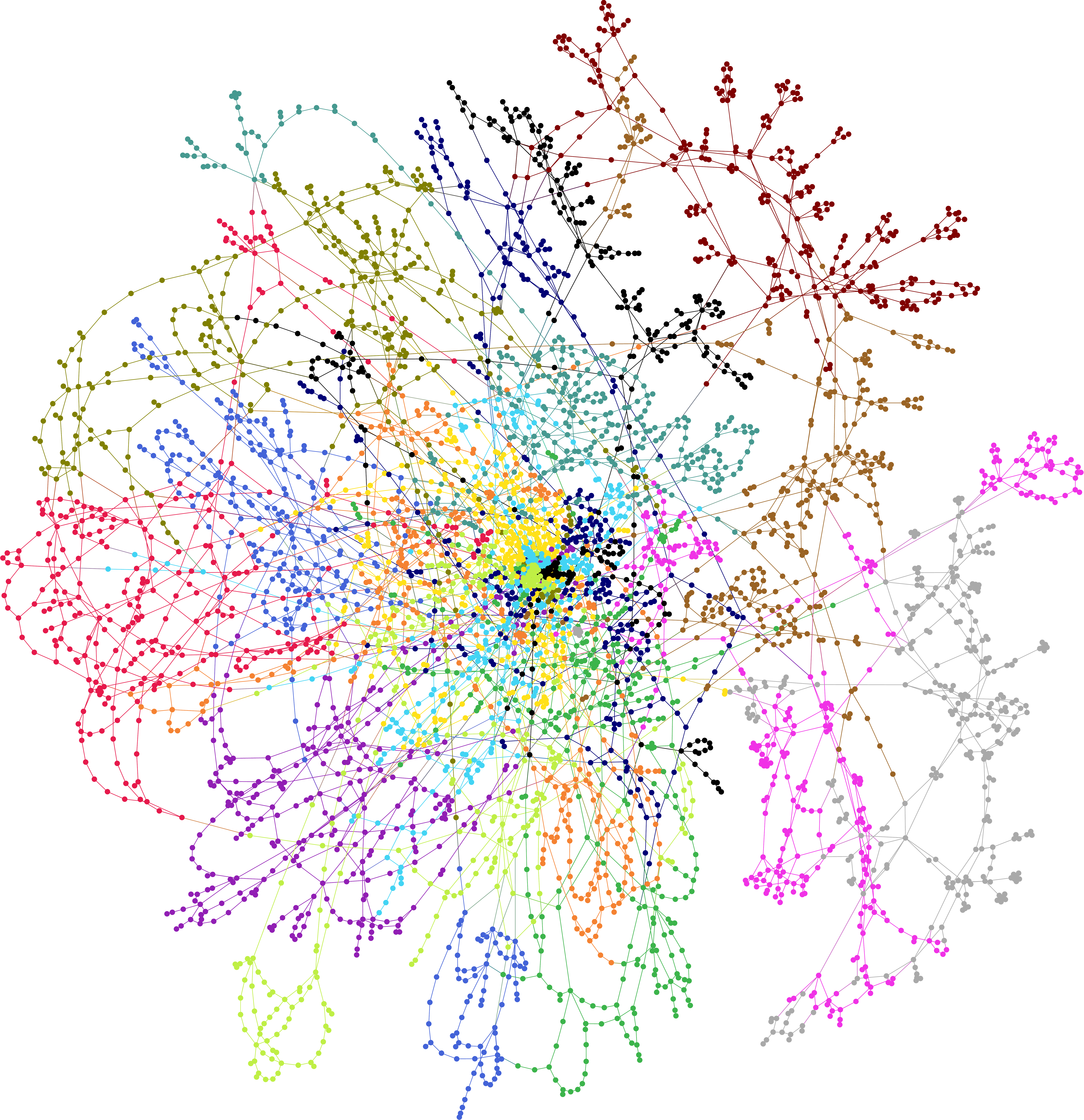}
    \includegraphics[width=.19\textwidth]{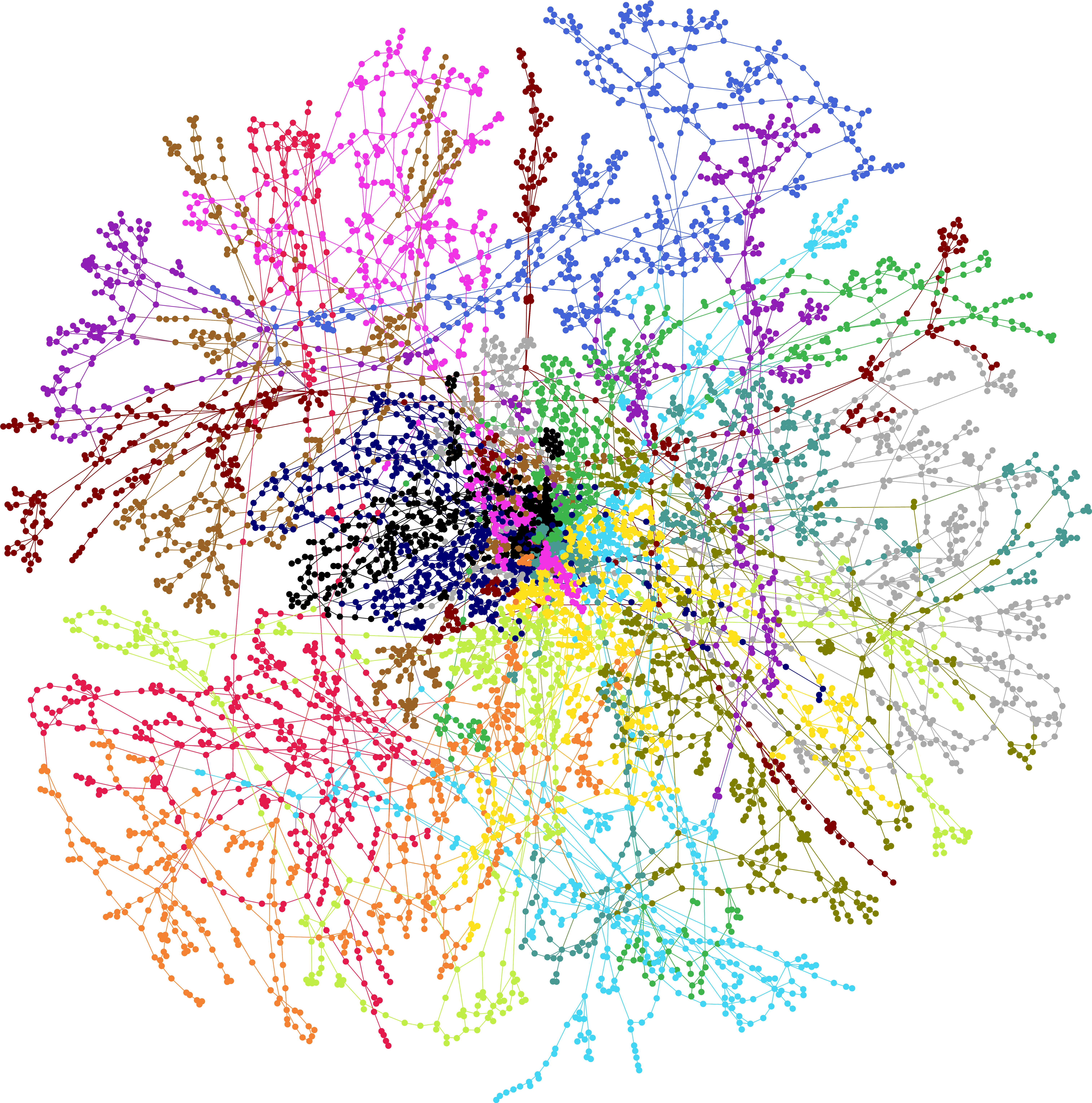}\\
    \includegraphics[width=.19\textwidth]{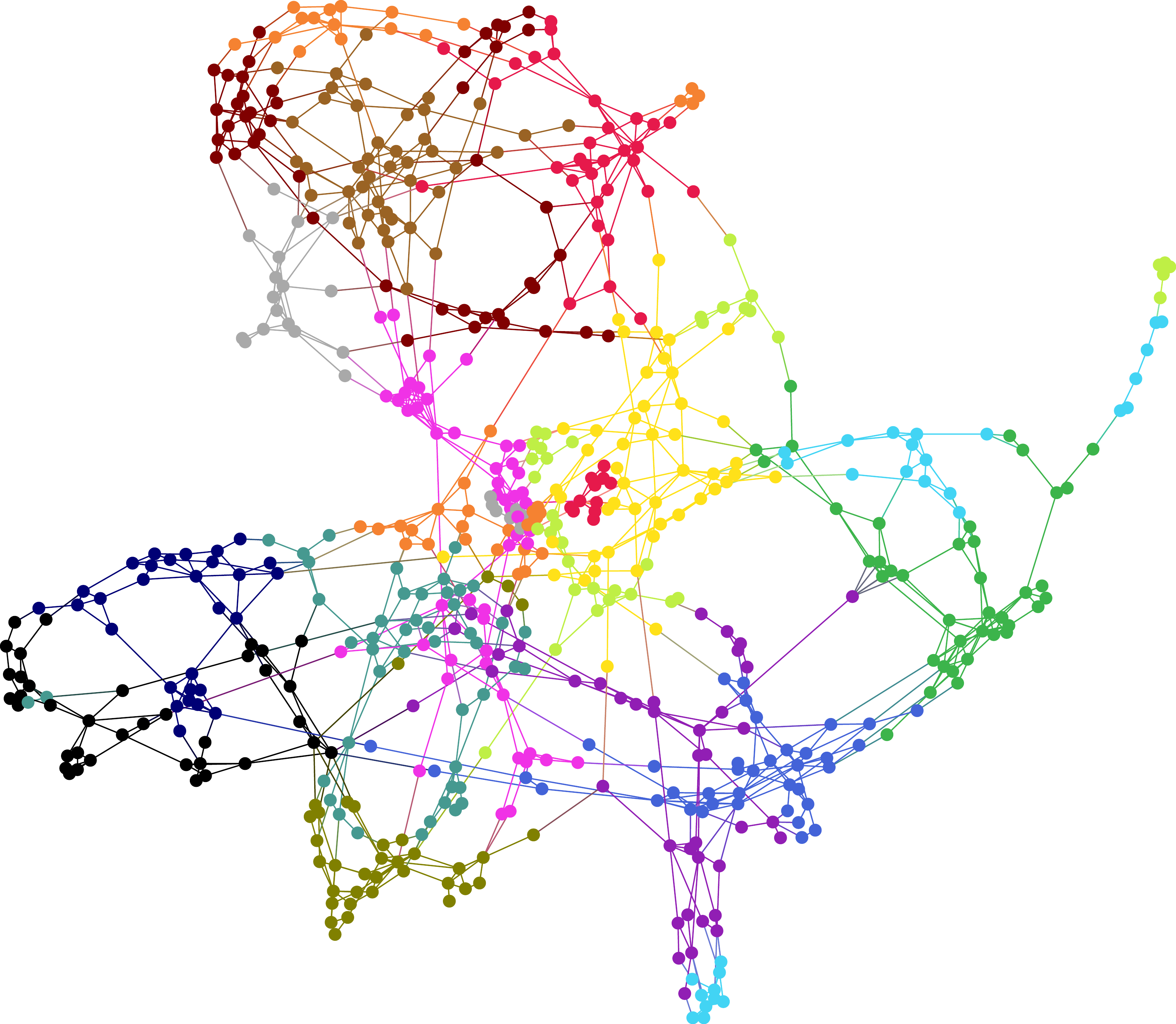}
    \includegraphics[width=.19\textwidth]{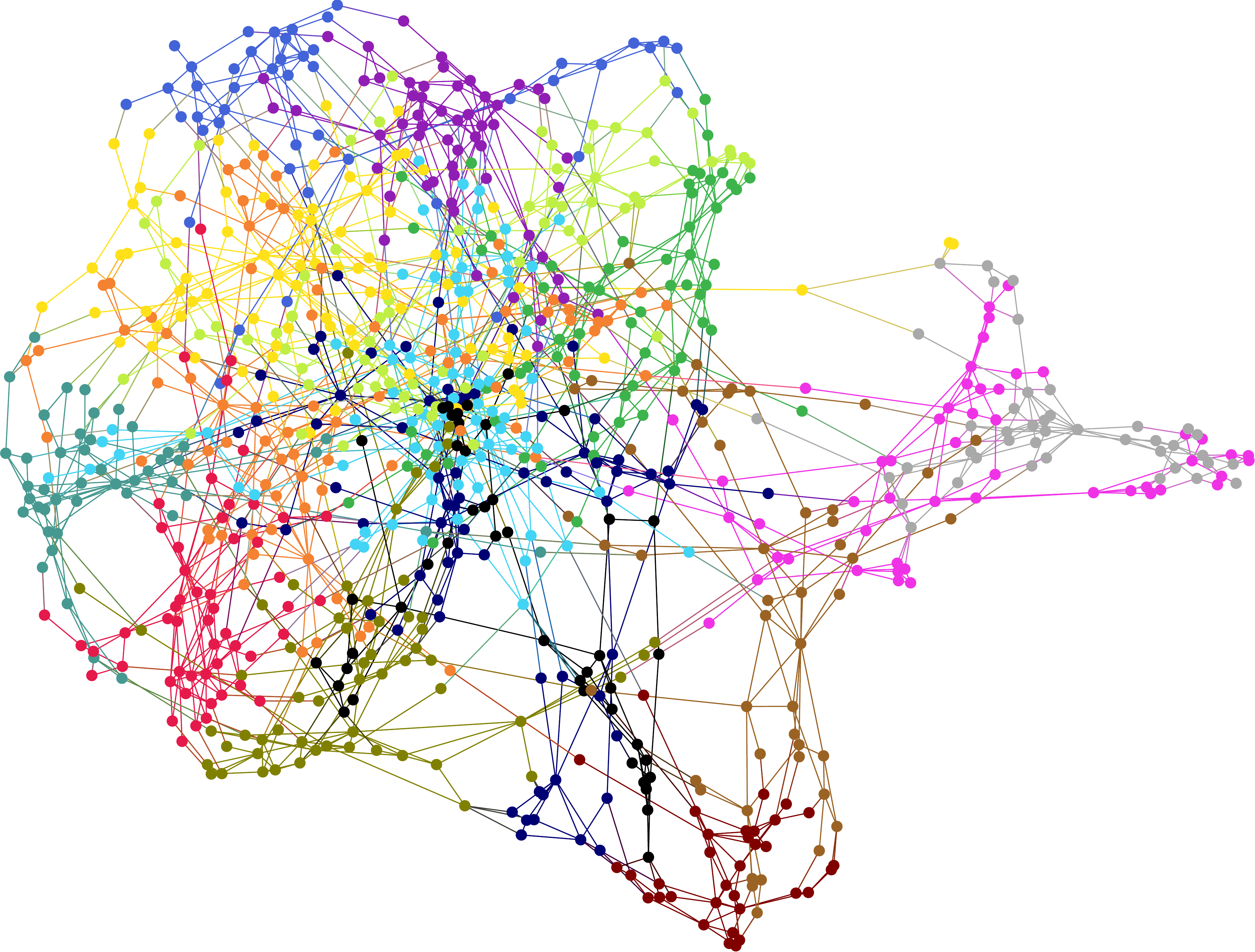}
    \includegraphics[width=.19\textwidth]{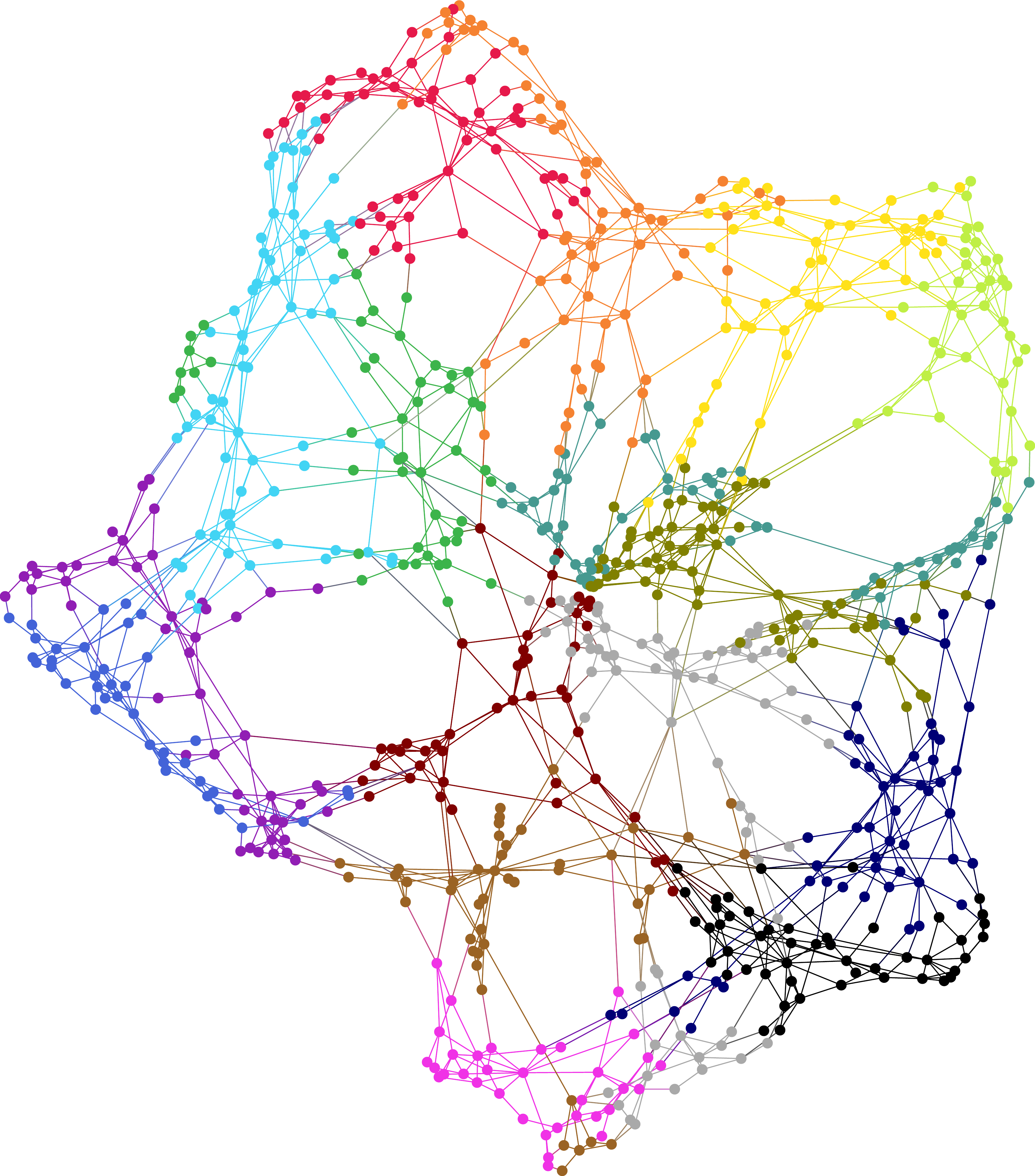}
    \caption{Left to right: Network A, Network C, Network E. Top to bottom: original network, coarse network.}
    \label{fig:networks}
  \end{center}
\end{figure*}

  We apply the proposed hierarchical scheme to solve \eqref{eqn:cs}. {\cm We use a set of network data selected from the benchmark library PGlib v18.08 \cite{pglib}. The characteristics of the networks are described in Table \ref{table:net} and the topologies for a subset of them are shown in Fig. \ref{fig:networks}.} To compare the performance, we also implement a centralized and a purely decentralized ADMM scheme (without supervisory layer). The original (centralized) problem \eqref{eqn:cs} and the subproblems are solved using the nonlinear programming solver {\tt Ipopt}\cite{Wachter2006}. The schemes are implemented in {\tt Julia} and leverage the algebraic modeling capabilities of {\tt JuMP} \cite{dunning2017jump}. We use {\tt METIS} for graph partitioning \cite{karypis1995metis}. The {\cm message passing interface (MPI) is used for parallel implementation. The case study is run on 17 processor cores (1 main process and 16 worker processes) using Intel(R) Xeon(R) CPU E5-2698 v3 @2.30GHz}.

\begin{table*}[!htbp]
  \caption{Network characteristics}
  \centering
  \label{table:net}
  \begin{tabular}{|c|c|c|c|c|c|c|c|c|c|c|c|}
    \hline
    &       & \multicolumn{5}{c|}{Full network}         & \multicolumn{5}{c|}{Coarse network} \\
    \hline
    Label & Name  & Nodes & Edges & Generators & Vars. & Constr. & Nodes & Edges & Generators & Var.  & Constr. \\
    \hline
    A     & case2853\_sdet & 2,853 & 3,635 & 946   & 20,490 & 22,528 & 865   & 1,512 & 946   & 7,158 & 6,732 \\
    \hline
    B     & case3120sp\_k & 3,120 & 3,684 & 505   & 21,158 & 24,068 & 1,100 & 1,695 & 505   & 5,462 & 5,616 \\
    \hline
    C     & case4661\_sdet & 4,661 & 5,751 & 1,176 & 33,212 & 37,107 & 1,697 & 2,750 & 1,176 & 10,778 & 10,662 \\
    \hline
    D     & case6468\_rte & 6,468 & 8,065 & 1,295 & 43,532 & 49,464 & 1,917 & 3,287 & 1,295 & 9,758 & 9,346 \\
    \hline
    E     & case6515\_rte & 6,515 & 8,104 & 1,388 & 43,928 & 49,682 & 1,874 & 3,248 & 1,388 & 11,708 & 11,638 \\
    \hline
  \end{tabular}%
\end{table*}

In order to formulate the coarse problems, one can use the aggregation scheme discussed in \eqref{eqn:coarse_functions}. Equivalently, one can consider the aggregation procedure as constructing the aggregated data for the coarse graph and formulating the coarse problem using the coarse graph and the aggregated data. By applying that aggregation scheme to \eqref{eqn:cs}, we find that the data should be aggregated as:
  \begingroup
  \allowdisplaybreaks
  \begin{subequations}\label{eqn:aggdata}
    \begin{align}
      G_{i^c j^c}^c&=
      \begin{cases}
        \sum_{i\in \tilde{\mathcal{V}}_{i^c} }G_{ii} + \sum_{\{i,j\}\in\mathcal{E}_{i^c,i^c}} 2G_{ij}  &\text{if }i^c=j^c\\
        \sum_{\{i,j\}\in \mathcal{E}_{i^c,j^c}} G_{ij}  &\text{if }i^c\neq j^c
      \end{cases}\\
      B_{i^c j^c}^c&=
      \begin{cases}
        \sum_{i\in \tilde{\mathcal{V}}_{i^c} }B_{ii} + \sum_{\{i,j\}\in\mathcal{E}_{i^c,i^c}} 2B_{ij}  &\text{if }i^c=j^c\\
        \sum_{\{i,j\}\in \mathcal{E}_{i^c,j^c}} B_{ij}  &\text{if }i^c\neq j^c
      \end{cases}\\
      P_L^c(i^c)&=\sum_{i\in \tilde{\mathcal{V}}_{i^c}} P_L(i),\quad
      Q_L^c(i^c)=\sum_{i\in\tilde{\mathcal{V}}_{i^c}} Q_L(i),
    \end{align}
  \end{subequations}
  \endgroup
  where $\mathcal{E}_{i^c,j^c}$ is defined in \eqref{eqn:Ejc}. See Appendix \ref{apx:agg} for a detailed derivation of \eqref{eqn:aggdata}. By formulating the OPF problem using the coarse graph $\mathcal{G}^c(\mathcal{V}^c,\mathcal{E}^c)$ and the aggregated data \eqref{eqn:aggdata}, one can construct the coarse problem. The characteristics of the coarse networks are described in Table \ref{table:net} and Fig. \ref{fig:networks}. Furthermore, the dimensions of the original and coarse problems (the number of variables and the number of equlaity and inequality constraints) are compared in Table \ref{table:net}.

We use the following metrics as stopping criteria:
\begin{align}
  \text{stop if}\quad\Vert r^\ell \Vert<\epsilon^{\text{pr}}\quad\text{and}\quad \Vert s^\ell \Vert<\epsilon^{\text{du}}.
\end{align}
The following formula suggested in \cite{Boyd2010} is used to set the primal tolerance $\epsilon^{\text{pr}}$ and dual tolerance $\epsilon^{\text{du}}$.
\begin{subequations}
  \begin{align}
    \epsilon^{\text{pr}}&=\sqrt{n_x}\epsilon^{\text{abs}}+\epsilon^{\text{rel}}\max\{\Vert Ax^\ell\Vert ,\Vert Bz^\ell\Vert\}\\
    \epsilon^{\text{du}}&=\sqrt{n_y}\epsilon^{\text{abs}}+\epsilon^{\text{rel}}\Vert A^T y^\ell \Vert
  \end{align}
\end{subequations}
where $\epsilon^{\text{abs}}$ and $\epsilon^{\text{rel}}$ are user-defined algorithm parameters. The algorithmic parameters used in the study are: $\rho=10^6$ to $10^7$ (depending on the instance) and $\epsilon^{\text{abs}}=\epsilon^{\text{rel}}=5\times10^{-4}$.


{\cm The results are summarized in Table \ref{table:results} and Fig. \ref{fig:results}. We compare the performance of three different solution schemes: centralized, purely decentralized, and hierarchical. Performance is measured in terms of objective value and solution time. For decentralized and hierarchical schemes, the progress of primal and dual residuals, objective values, and augmented Lagrangian values are compared.

\begin{table*}[!htbp]
  \caption{Performance of decentralized and hierarchical schemes.}
  \centering
  \label{table:results}
  \begin{tabular}{|c|c|c|c|c|c|c|c|c|c|}
      \hline
      & \multicolumn{2}{c|}{Centralized} & \multicolumn{3}{c|}{Purely Decentralized} & \multicolumn{4}{c|}{Hierarchical} \\
      \hline
      Label& Time (s) & Obj. ($\times 10^6$) & Time (s) & ADMM step & Obj. ($\times 10^6$)& Coarse Time (s) & Overall Time (s) & ADMM step & Obj. ($\times 10^6$) \\
      \hline
A & 5.9564 & 2.5530 & 40.8169 & 46 & 2.5619 & 1.3114 & 27.9388 & 31 & 2.5541\\\hline
B & 3.2632 & 1.6464 & 157.4092 & 157 & 1.6276 & 2.4252 & 45.0053 & 43 & 1.6438\\\hline
C & 15.8135 & 2.7516 & 229.7916 & 109 & 2.9108 & 4.2938 & 88.0437 & 40 & 2.7779\\\hline
D & 20.5671 & 1.7595 & 242.3939 & 61 & 1.7909 & 2.4635 & 152.2255 & 38 & 1.7979\\\hline
E & 22.5522 & 2.5052 & 173.7644 & 43 & 2.9728 & 2.5857 & 109.7643 & 26 & 2.7837\\\hline
    \end{tabular}%
\end{table*}
\begin{figure*}[!htbp]
  \begin{center}
    \includegraphics[width=\textwidth]{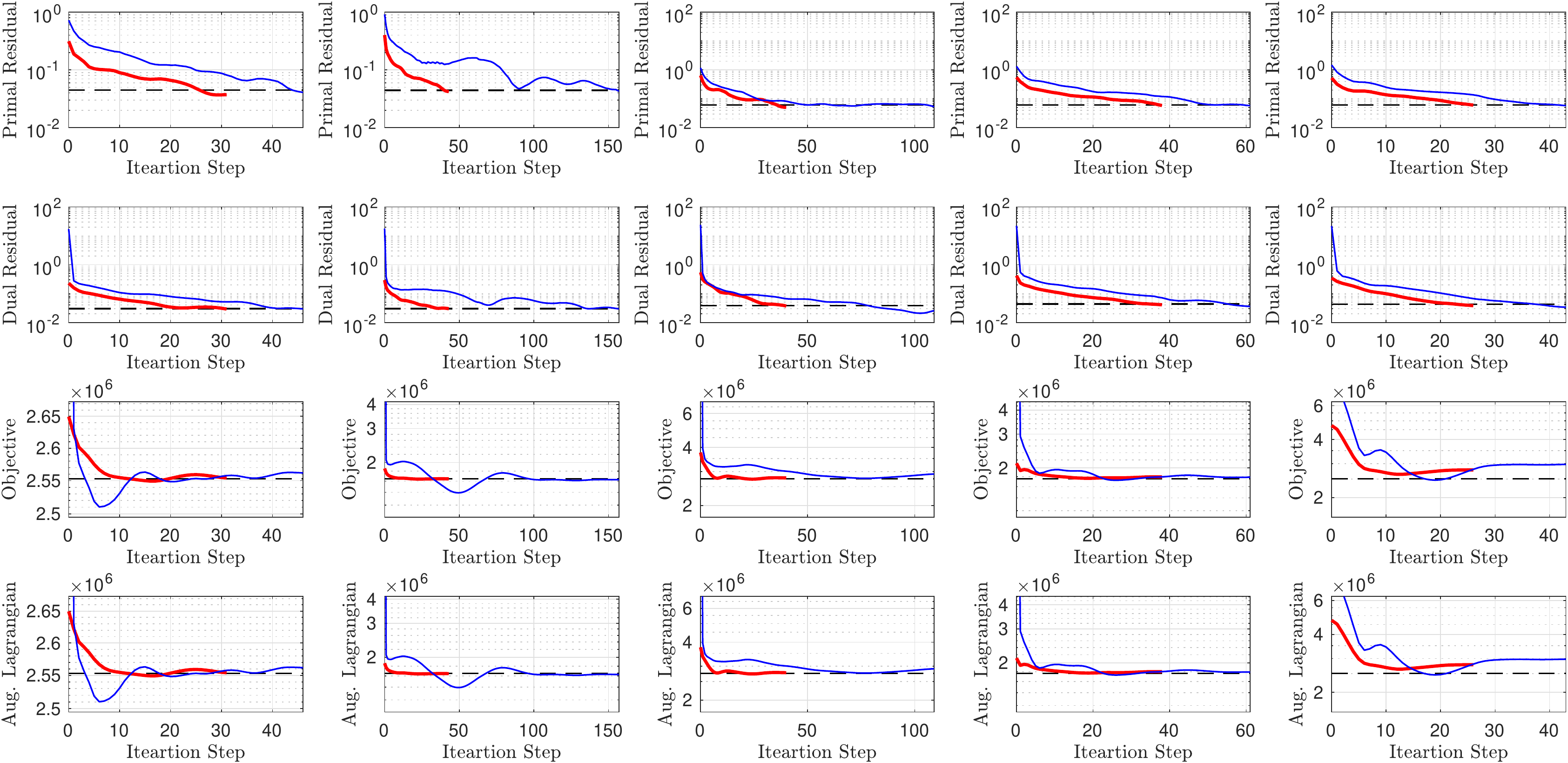}
\caption{Blue lines are trajectories for the decentralized ADMM scheme and red lines are trajectories for the hierarchical scheme. From left to right: Network A, Network B, Network C, Network D, Network E. From top to bottom: primal residual, and dual residual, objective value, and augmented Lagrangian.}
\label{fig:results}
\end{center}
\end{figure*}

We have found that, for the test instances, both the decentralized scheme and the hierarchical scheme converge to a stationary point if a large enough penalty parameter $\rho$ is chosen. We have also found relatively small gaps (on the order of 1\%) between the objective values of the centralized, decentralized, and hierarchical schemes.}

By comparing the evolution of primal and dual residuals, we can clearly see the benefit of the hierarchical architecture over the decentralized one. In particular, the convergence trajectories are smoother and the number of coordination steps are drastically reduced (by 18-74\%). Furthermore, we can observe that the computation times for coarse problems are significantly shorter than the computation times for ADMM coordination. Accordingly, as can be seen from the results in Table \ref{table:results}, the overall computation times for the hierarchical scheme are always smaller than the computation times for the purely decentralized scheme. This indicates that the hierarchical scheme can achieve a faster coordination among the subnetworks. 
  
{\cm We found that, for the instances considered, the hierarchical scheme is not as fast as the centralized scheme. Benefits in computing time are expected in larger instances (e.g., when linear algebra operations in a centralized solver reach capacity). Consequently, we only present the centralized solution in order to validate performance. We also highlight that benefits of decentralization go beyond solution time (e.g., information privacy, reduced communication requirements, and flexibility in implementation).}

\section{Conclusions}\label{sec:conclusion}
We have presented a hierarchical optimization architecture for large-scale power networks that overcomes limitations of fully centralized and fully decentralized architectures. The top layer of the architecture uses a coarse representation of the entire network while the bottom layer is composed of a family of decentralized optimization agents each operating on a network subdomain at full resolution. We show that state and dual information obtained from the top layer can be used to accelerate the coordination of the decentralized optimization agents and to recover optimality for the entire system. We use an alternating direction method of multipliers (ADMM) framework to drive coordination of the decentralized agents. We provide procedures to construct coarse representations for power networks and demonstrate that the hierarchical architecture can handle large systems. 

\appendix

\subsection{Derivation of Stationarity Conditions}\label{apx:KKT}  
Let $(x^{\ell+1},\lambda^{\ell+1},\nu^{\ell+1})$ be a primal-dual solution of \eqref{eqn:ADMM_x} at $\ell$th iteration such that satisfies some constraint qualification.
We have that $(x^{\ell+1},\lambda^{\ell+1},\nu^{\ell+1})$ satisfies the stationarity of \eqref{eqn:ADMM_x}.
Thus, we can establish \eqref{eqn:pr_fs_1}, \eqref{eqn:pr_fs_2}, \eqref{eqn:du_fs}, and \eqref{eqn:cmpl_slk}. Using the first order conditions of \eqref{eqn:ADMM_x}-\eqref{eqn:ADMM_z} and \eqref{eqn:ADMM_y} we obtain:
\begin{subequations}\label{eqn:ADMM_KKT}
  \begin{align}
    0= &\nabla_x f(x^{\ell+1})+A^T y^{\ell}+\rho A^T (Ax^{\ell+1}+Bz^{\ell})\\\nonumber
&\quad+\nabla_x g(x^{\ell+1})\lambda^{\ell+1}+\nabla_x h(x^{\ell+1})\nu^{\ell+1}\\
=&\nabla_x f(x^{\ell+1})+A^T y^{\ell+1}+\rho A^T B(-z^{\ell+1}+z^{\ell})\\\nonumber
&\quad+\nabla_x g(x^{\ell+1})\lambda^{\ell+1}+\nabla_x h(x^{\ell+1})\nu^{\ell+1}\\
0=& B^T y^{\ell} +\rho B^T (Ax^{\ell+1}+Bz^{\ell+1})=B^T y^{\ell+1}.\label{eqn:By=0}
\end{align}
\end{subequations}
Using \eqref{eqn:ADMM_KKT}, we can derive the following:
\begin{subequations}
\begin{align}
&\nabla_{x}\mathcal{L} (x^{\ell+1},z^{\ell+1},y^{\ell+1},\lambda^{\ell+1},\nu^{\ell+1})\\\nonumber
&=\nabla_x f(x^{\ell+1})+A^T y^{\ell+1}\\\nonumber
&\qquad+\nabla_x g(x^{\ell+1})\lambda^{\ell+1}+\nabla_x h(x^{\ell+1})\nu^{\ell+1}\\\nonumber
&=\rho A^T B(z^{\ell+1}-z^{\ell})=s^{\ell+1}
\\
&\nabla_{z}\mathcal{L} (x^{\ell+1},z^{\ell+1},y^{\ell+1},\lambda^{\ell+1},\nu^{\ell+1})=B^T y^{\ell+1}=0
\\
&\nabla_{y}\mathcal{L} (x^{\ell+1},z^{\ell+1},y^{\ell+1},\lambda^{\ell+1},\nu^{\ell+1})\\\nonumber
&=Ax^{\ell+1}+By^{\ell+1}=r^{\ell+1}.
\end{align}
\end{subequations}
This establishes conditions \eqref{eqn:stationary_1}, \eqref{eqn:stationary_2}, \eqref{eqn:pr_fs_3}.

\subsection{Coarse OPF Problem Fromulation}\label{apx:agg}
In this section, we derive the coarse OPF problem for \eqref{eqn:cs} using the coarse graph $\mathcal{G}^c (\mathcal{V}^c, \mathcal{E}^c)$ and the subpartition structure $\{\tilde{\mathcal{V}}_1,\cdots,\tilde{\mathcal{V}}_{K^c}\}$. For each $i^c=1,2,\cdots,K^c$, we treat the voltage angles and amplitudes of all the nodes in $\tilde{\mathcal{V}}_{i^c}$ as the same. In particular, we enforce that:
  \begin{align}
    \theta^c(i^c) &= \theta(i),\quad  V^c(i^c) = V(i),\quad\forall i\in \tilde{V}_{i^c},
  \end{align}
  where the voltage angle and the voltage amplitude of node $i^c\in\mathcal{V}^c$ are denoted by $V^c(i^c)$ and $\theta^c(i^c)$, respectively.

  We augment the power generation vectors in each subpartition. In particular, we first take the union of the generator sets to form $\mathcal{W}^c_{i^c} := \bigcup_{i\in \tilde{V}_{i^c}} \mathcal{W}_i$, and consider
      \begin{align}
        P(\mathcal{W}^c_{i^c}) &= \{P(j)\}_{j\in \mathcal{W}^c_{i^c}},\quad 
        Q(\mathcal{W}^c_{i^c}) = \{Q(j)\}_{j\in \mathcal{W}^c_{i^c}}
      \end{align}
    as the power generation vectors for $i^c\in\mathcal{V}^c$. The state vector for each node can be defined as follows:
    \begin{align}
      x^c_{i^c}:=(\theta^c_{i^c},V^c_{i^c},P(\mathcal{W}^c_{i^c}),Q(\mathcal{W}^c_{i^c})).
    \end{align}
    
 {\em Objective function:}  By \eqref{eqn:crsf}, the objective function in subpartition $\tilde{\mathcal{V}}^c_{i^c}$ is aggregated as:
  \begin{align}\label{eqn:aggobj}
    f_{c,(i^c)}(x^c(i^c))=\sum_{i\in \tilde{\mathcal{V}}_{i^c} }\left[\sum_{j\in \mathcal{W}_{i}} c_{j} P(j)\right] =     \sum_{j\in \mathcal{W}^c_{i^c}} c_{j} P(j)
  \end{align}
  With the summation over $i^c\in\mathcal{V}^c$, one can construct the overall objective function. We can see that the overall objective function does not change by the graph aggregation.

  {\em Equality constraints:} By \eqref{eqn:crsg}, the active power flow constraints \eqref{eqn:cs_b} in subpartition $\tilde{\mathcal{V}}^c_{i^c}$ are aggregated as follows:
    \begin{align}
      &\sum_{i\in\tilde{\mathcal{V}}_{i^c}} \left[P_L(i)+\sum_{j\in \mathcal{W}_i} P(j)\right]\\
      \nonumber &=\sum_{i\in\tilde{\mathcal{V}}_{i^c}} \Bigg[ \sum_{j\in \mathcal{N}[i]}G_{ij}V (i)V (j)\cos\left(\theta (i)-\theta (j)\right)\\
        \nonumber &+\sum_{j\in \mathcal{N}[i]}B_{ij}V (i)V (j)\sin\left(\theta (i)-\theta (j)\right)\Bigg],\;\forall i^c\in\mathcal{V}^c
    \end{align}
    By inspecting the algebraic structure, we can write:
    {\footnotesize
    \begin{align}
      &\left[\sum_{i\in\tilde{\mathcal{V}}_{i^c}} P_L(i)\right]+\sum_{j\in \mathcal{W}^c_{i^c}} P(j)\\
      \nonumber &=\left[\sum_{i\in \tilde{\mathcal{V}}_{i^c} }G_{ii} +  \sum_{\{i,j\}\in\mathcal{E}_{i^c,i^c}}2G_{ij}\right] V^c(i^c)^2\\
      \nonumber &+\sum_{j^c\in \mathcal{N}_{\mathcal{G}^c}(i^c)}
      \left[\sum_{\{i,j\}\in\mathcal{E}_{i^c,j^c}} G_{ij}\right]
      V^c (i^c)V^c (j^c)\cos\left(\theta^c (i^c)-\theta^c (j^c)\right)\\
      \nonumber &+\sum_{j^c\in \mathcal{N}_{\mathcal{G}^c}(i^c)}
      \left[\sum_{\{i,j\}\in\mathcal{E}_{i^c,j^c}} B_{ij}\right]
      V^c (i^c)V^c (j^c)\sin\left(\theta^c (i^c)-\theta^c (j^c)\right)
    \end{align}}
    Similarly, by aggregating the reactive power flow constraints \eqref{eqn:cs_c} in subpartition $\tilde{\mathcal{V}}^c_{i^c}$, we obatin:
    {\footnotesize
      \begin{align}
        &\left[\sum_{i\in\tilde{\mathcal{V}}_{i^c}} Q_L(i)\right]+\sum_{j\in \mathcal{W}^c_{i^c}} Q(j)\\
        \nonumber &=\sum_{j^c\in \mathcal{N}_{\mathcal{G}^c}(i^c)}
        \left[\sum_{\{i,j\}\in\mathcal{E}_{i^c,j^c}} G_{ij}\right]
        V^c (i^c)V^c (j^c)\sin\left(\theta^c (i^c)-\theta^c (j^c)\right)\\
        \nonumber &-\left[\sum_{i\in \tilde{\mathcal{V}}_{i^c} }B_{ii} +  \sum_{\{i,j\}\in\mathcal{E}_{i^c,i^c}} 2B_{ij}\right] V^c(i^c)^2\\
        \nonumber &-\sum_{j^c\in \mathcal{N}_{\mathcal{G}^c}(i^c)}
        \left[\sum_{\{i,j\}\in\mathcal{E}_{i^c,j^c}} B_{ij}\right]
        V^c(i^c)V^c(j^c)\cos\left(\theta^c (i^c)-\theta^c (j^c)\right)
    \end{align}}
  We observe that defining $G^c$, $B^c$, $P^c_L$, and $Q^c_L$ as \eqref{eqn:aggdata} yields:
  \begin{subequations}
    \begin{align}
      &P^c_L(i^c)+\sum_{j\in \mathcal{W}^c_{i^c}} P(j)=\\
      &\sum_{j^c\in \mathcal{N}_{\mathcal{G}^c}[i^c]}G^c_{i^cj^c}V^c (i^c)V^c (j^c)\cos\left(\theta^c (i^c)-\theta^c (j^c)\right)\nonumber\\
      &+\sum_{j^c\in\mathcal{N}_{\mathcal{G}^c}[i^c]}B^c_{i^cj^c}V^c (i^c)V^c (j^c)\sin\left(\theta^c (i^c)-\theta^c (j^c)\right)\nonumber\\
      &Q^c_L(i^c)+\sum_{j\in \mathcal{W}^c_{i^c}} Q(j)=\\
      &\sum_{j^c \in \mathcal{N}_{\mathcal{G}^c}[i^c]}G^c_{i^c j^c}V^c (i^c)V^c (j^c)\sin\left(\theta^c (i^c)\ - \theta^c (j^c)\right)\nonumber\\
      &-\sum_{j^c\in\mathcal{N}_{\mathcal{G}^c}[i^c]}B^c_{i^c j^c}V^c (i^c)V^c (j^c)\cos\left(\theta^c (i^c)\ - \theta^c (j^c)\right).\nonumber
    \end{align}
  \end{subequations}
  Thus, we recover the power flow equation form for the coarse graph $\mathcal{G}^c(\mathcal{V}^c,\mathcal{E}^c)$.

  The equality constraints for reference nodes \eqref{eqn:cs_h} can be simply reduced to:
  \begin{align}
    \theta^c (i^c) = 0
  \end{align}
  for $i^c\in\mathcal{V}^{c,\text{ref}}$, where $\mathcal{V}^{c,\text{ref}}:= \{i^c\in\mathcal{V}^c\mid \tilde{\mathcal{V}}_{i^c}\cap\mathcal{V}^{\text{ref}}\neq \emptyset\}$.

  {\em Inequality constraints:} Constraints \eqref{eqn:cs_d}-\eqref{eqn:cs_e} does not change with the aggregation because the aggregation does not affect the generators.
  By aggregation scheme \eqref{eqn:crsh-e}, one can obtain:
  \begin{align}
    \theta^{\min} \leq \theta^c(i^c)-\theta^c(j^c) \leq \theta^{\max},\quad \{i^c,j^c\}\in\mathcal{E}^c.
  \end{align}
  Lastly, with the \eqref{eqn:crsh-v}, \eqref{eqn:cs_g} reduces to:
  \begin{align}\label{eqn:aggvoltage}
    V^{\min} \leq V^c(i^c)\leq V^{\max}, \quad i^c\in\mathcal{V}^c.
  \end{align}
  
By combining \eqref{eqn:aggobj}-\eqref{eqn:aggvoltage}, we find that formulating the OPF problem with aggregated graph $\mathcal{G}^c(\mathcal{V}^c,\mathcal{E}^c)$ and the aggregated data \eqref{eqn:aggdata} is equivalent to aggregating the objective function and constraints based on the summation procedure \eqref{eqn:coarse_functions}. By applying the lifting procedure, one can obtain the form of \eqref{eqn:coarse}.

\ifCLASSOPTIONcaptionsoff
  \newpage
\fi

\bibliographystyle{IEEEtran}
\bibliography{multigrid}

\vspace{-0.4in}
\begin{IEEEbiographynophoto}{Sungho Shin}
received his B.S. in chemical engineering and mathematics from Seoul National University, South Korea in 2016. He is currently a Ph.D. candidate in the Department of Chemical and Biological Engineering at the University of Wisconsin-Madison. His research interests include control theory, optimization algorithms, and complex networks.
\end{IEEEbiographynophoto}

\vspace{-0.4in}
\begin{IEEEbiographynophoto}{Philip Hart}
graduated from Clarkson University, NY, USA in 2011, with a BSc in Electrical Engineering and a minor in Sustainable Energy Systems Engineering. He received a MSrc degree in Electrical Engineering at the University of Wisconsin-Madison, WI, USA in 2013. He is currently a PhD candidate in the Electrical and Computer Engineering Department at the University of Wisconsin-Madison, pursuing research interests in distributed generation, nonlinear microgrid dynamics, and model order reduction.
\end{IEEEbiographynophoto}

\vspace{-0.3in}
\begin{IEEEbiographynophoto}{Thomas Jahns} (S'73-M'79-SM'91-F'93) received the S.B., S.M., and Ph.D. degrees in electrical engineering from the Massachusetts Institute of Technology, Cambridge, MA (USA) in 1974 and 1978. In 1998, he joined the Department of Electrical and Computer Engineering, University of Wisconsin-Madison, as a Grainger Professor of Power Electronics and Electric Machines, where he is currently a Co-Director of the Wisconsin Electric Machines and Power Electronics Consortium (WEMPEC).  Prior to joining UW, he worked at GE Corporate Research and Development (now GE Global Research Center), Niskayuna, NY, for 15 years. His current research interests include high-performance permanent-magnet synchronous machines, electric traction drives, and distributed energy resources, including microgrids.  Dr. Jahns received the 2005 IEEE Nikola Tesla Technical Field Award and the IAS Outstanding Achievement Award in 2011.  He was elected as a member of the National Academy of Engineering in 2015.
\end{IEEEbiographynophoto}

\vspace{-0.3in}
\begin{IEEEbiographynophoto}{Victor M. Zavala}
is the Baldovin-DaPra Associate Professor in the Department of Chemical and Biological Engineering at the University of Wisconsin-Madison. He holds a B.Sc. degree from Universidad Iberoamericana and a Ph.D. degree from Carnegie Mellon University, both in chemical engineering. He is an associate editor for the Journal of Process Control and a technical editor of Mathematical Programming Computation. His research interests are in the areas of energy systems, high-performance computing, stochastic optimization, and predictive control.\end{IEEEbiographynophoto}
\end{document}